\input amstex
\input amsppt.sty

\magnification1200
\hsize13cm
\vsize19cm

\TagsOnRight

\def\StanAP{9}
\def\ReZeAA{8}
\def\ReVeAA{7}
\def\RegeAE{6}
\def\RegeAD{5}
\def\PoReAA{4}
\def\MacdAC{3}
\def\JanSAA{2}
\def\BessAC{1}

\def\HP{H\!P}
\def\SR{S\kern-1ptR}
\def\SQ{S\kern-.6ptQ}
\def\leaderfill{\leaders\hbox to 1em{\hss.\hss}\hfill}

\def\LL{\leavevmode\setbox0=\hbox{L}\hbox to\wd0{\hss\char'40L}}

\def\la{\lambda}

            %used for crossreferencing, Tex should ignore.
             %used for refencing (section-numbers)
          %used for new-section numbers

\def\today{\ifcase\month\or
 January\or February\or March\or April\or May\or June\or
 July\or August\or September\or October\or November\or December\fi
 \space\number\day, \number\year}
 %zum Nummerieren

\def\({\left(}
\def\){\right)}
\def\[{\left[}
\def\]{\right]}

\def\3{\ss}
\catcode`\@=11
\def\dddot#1{\vbox{\ialign{##\crcr
      .\hskip-.5pt.\hskip-.5pt.\crcr\noalign{\kern1.5\p@\nointerlineskip}
      $\hfil\displaystyle{#1}\hfil$\crcr}}}

\newif\iftab@\tab@false
\newif\ifvtab@\vtab@false
\def\tab{\bgroup\tab@true\vtab@false\vst@bfalse\Strich@false%
   \def\\{\global\hline@@false%
     \ifhline@\global\hline@false\global\hline@@true\fi\cr}
   \edef\l@{\the\leftskip}\ialign\bgroup\hskip\l@##\hfil&&##\hfil\cr}
\def\endtab{\cr\egroup\egroup}
\def\vtab{\vtop\bgroup\vst@bfalse\vtab@true\tab@true\Strich@false%
   \bgroup\def\\{\cr}\ialign\bgroup&##\hfil\cr}
\def\endvtab{\cr\egroup\egroup\egroup}
\def\stab{\D@cke0.5pt\null 
 \bgroup\tab@true\vtab@false\vst@bfalse\Strich@true\Let@@\vspace@
 \normalbaselines\offinterlineskip
  \openup\spreadmlines@
 \edef\l@{\the\leftskip}\ialign
 \bgroup\hskip\l@##\hfil&&##\hfil\crcr}
\def\endstab{\crcr\egroup
 \egroup}
\newif\ifvst@b\vst@bfalse
\def\vstab{\D@cke0.5pt\null
 \vtop\bgroup\tab@true\vtab@false\vst@btrue\Strich@true\bgroup\Let@@\vspace@
 \normalbaselines\offinterlineskip
  \openup\spreadmlines@\bgroup}
\def\endvstab{\crcr\egroup\egroup
 \egroup\tab@false\Strich@false}

\newdimen\htstrut@
\htstrut@8.5\p@
\newdimen\htStrut@
\htStrut@12\p@
\newdimen\dpstrut@
\dpstrut@3.5\p@
\newdimen\dpStrut@
\dpStrut@3.5\p@
\def\openup{\afterassignment\@penup\dimen@=}
\def\@penup{\advance\lineskip\dimen@
  \advance\baselineskip\dimen@
  \advance\lineskiplimit\dimen@
  \divide\dimen@ by2
  \advance\htstrut@\dimen@
  \advance\htStrut@\dimen@
  \advance\dpstrut@\dimen@
  \advance\dpStrut@\dimen@}
\def\Let@@{\relax\iffalse{\fi%
    \def\\{\global\hline@@false%
     \ifhline@\global\hline@false\global\hline@@true\fi\cr}%
    \iffalse}\fi}
\def\matrix{\null\,\vcenter\bgroup
 \tab@false\vtab@false\vst@bfalse\Strich@false\Let@@\vspace@
 \normalbaselines\openup\spreadmlines@\ialign
 \bgroup\hfil$\m@th##$\hfil&&\quad\hfil$\m@th##$\hfil\crcr
 \Mathstrut@\crcr\noalign{\kern-\baselineskip}}
\def\endmatrix{\crcr\Mathstrut@\crcr\noalign{\kern-\baselineskip}\egroup
 \egroup\,}
\def\smatrix{\D@cke0.5pt\null\,
 \vcenter\bgroup\tab@false\vtab@false\vst@bfalse\Strich@true\Let@@\vspace@
 \normalbaselines\offinterlineskip
  \openup\spreadmlines@\ialign
 \bgroup\hfil$\m@th##$\hfil&&\quad\hfil$\m@th##$\hfil\crcr}
\def\endsmatrix{\crcr\egroup
 \egroup\,\Strich@false}
\newdimen\D@cke
\def\Dicke#1{\global\D@cke#1}
\newtoks\tabs@\tabs@{&}
\newif\ifStrich@\Strich@false
\newif\iff@rst

\def\Stricherr@{\iftab@\ifvtab@\errmessage{\noexpand\s not allowed
     here. Use \noexpand\vstab!}%
  \else\errmessage{\noexpand\s not allowed here. Use \noexpand\stab!}%
  \fi\else\errmessage{\noexpand\s not allowed
     here. Use \noexpand\smatrix!}\fi}
\def\format{\ifvst@b\else\crcr\fi\egroup\iffalse{\fi\ifnum`}=0 \fi\format@}
\def\format@#1\\{\def\preamble@{#1}%
 \def\Str@chfehlt##1{\ifx##1\s\Stricherr@\fi\ifx##1\\\let\Next\relax%
   \else\let\Next\Str@chfehlt\fi\Next}%
 \def\c{\hfil\noexpand\ifhline@@\hbox{\vrule height\htStrut@%
   depth\dpstrut@ width\z@}\noexpand\fi%
   \ifStrich@\hbox{\vrule height\htstrut@ depth\dpstrut@ width\z@}%
   \fi\iftab@\else$\m@th\fi\the\hashtoks@\iftab@\else$\fi\hfil}%
 \def\r{\hfil\noexpand\ifhline@@\hbox{\vrule height\htStrut@%
   depth\dpstrut@ width\z@}\noexpand\fi%
   \ifStrich@\hbox{\vrule height\htstrut@ depth\dpstrut@ width\z@}%
   \fi\iftab@\else$\m@th\fi\the\hashtoks@\iftab@\else$\fi}%
 \def\l{\noexpand\ifhline@@\hbox{\vrule height\htStrut@%
   depth\dpstrut@ width\z@}\noexpand\fi%
   \ifStrich@\hbox{\vrule height\htstrut@ depth\dpstrut@ width\z@}%
   \fi\iftab@\else$\m@th\fi\the\hashtoks@\iftab@\else$\fi\hfil}%
 \def\s{\ifStrich@\ \the\tabs@\vrule width\D@cke\the\hashtoks@%
          \fi\the\tabs@\ }%
 \def\sa{\ifStrich@\vrule width\D@cke\the\hashtoks@%
            \the\tabs@\ %
            \fi}%
 \def\se{\ifStrich@\ \the\tabs@\vrule width\D@cke\the\hashtoks@\fi}%
 \def\cd{\hfil\noexpand\ifhline@@\hbox{\vrule height\htStrut@%
   depth\dpstrut@ width\z@}\noexpand\fi%
   \ifStrich@\hbox{\vrule height\htstrut@ depth\dpstrut@ width\z@}%
   \fi$\dsize\m@th\the\hashtoks@$\hfil}%
 \def\rd{\hfil\noexpand\ifhline@@\hbox{\vrule height\htStrut@%
   depth\dpstrut@ width\z@}\noexpand\fi%
   \ifStrich@\hbox{\vrule height\htstrut@ depth\dpstrut@ width\z@}%
   \fi$\dsize\m@th\the\hashtoks@$}%
 \def\ld{\noexpand\ifhline@@\hbox{\vrule height\htStrut@%
   depth\dpstrut@ width\z@}\noexpand\fi%
   \ifStrich@\hbox{\vrule height\htstrut@ depth\dpstrut@ width\z@}%
   \fi$\dsize\m@th\the\hashtoks@$\hfil}%
 \ifStrich@\else\Str@chfehlt#1\\\fi%
 \setbox\z@\hbox{\xdef\Preamble@{\preamble@}}\ifnum`{=0 \fi\iffalse}\fi
 \ialign\bgroup\span\Preamble@\crcr}
\newif\ifhline@\hline@false
\newif\ifhline@@\hline@@false
\def\hlinefor#1{\multispan@{\strip@#1 }\leaders\hrule height\D@cke\hfill%
    \global\hline@true\ignorespaces}
\def\Item "#1"{\par\noindent\hangindent2\parindent%
  \hangafter1\setbox0\hbox{\rm#1\enspace}\ifdim\wd0>2\parindent%
  \box0\else\hbox to 2\parindent{\rm#1\hfil}\fi\ignorespaces}
\def\ITEM #1"#2"{\par\noindent\hangafter1\hangindent#1%
  \setbox0\hbox{\rm#2\enspace}\ifdim\wd0>#1%
  \box0\else\hbox to 0pt{\rm#2\hss}\hskip#1\fi\ignorespaces}
\def\item"#1"{\par\noindent\hang%
  \setbox0=\hbox{\rm#1\enspace}\ifdim\wd0>\the\parindent%
  \box0\else\hbox to \parindent{\rm#1\hfil}\enspace\fi\ignorespaces}
\let\plainitem@\item
\catcode`\@=13

\catcode`\@=11
\font\tenln    = line10
\font\tenlnw   = linew10

\newskip\Einheit \Einheit=0.5cm
\newcount\xcoord \newcount\ycoord
\newdimen\xdim \newdimen\ydim \newdimen\PfadD@cke \newdimen\Pfadd@cke

%%%%%%%%%%%%%%%%%%%%%%%%%%%%%%%%%%%%%%%%%%%%%%%%%
%LaTeX counters, dimensions, variables for lines%
%%%%%%%%%%%%%%%%%%%%%%%%%%%%%%%%%%%%%%%%%%%%%%%%%
\newcount\@tempcnta
\newcount\@tempcntb

\newdimen\@tempdima
\newdimen\@tempdimb

\newdimen\@wholewidth
\newdimen\@halfwidth

\newcount\@xarg
\newcount\@yarg
\newcount\@yyarg
\newbox\@linechar
\newbox\@tempboxa
\newdimen\@linelen
\newdimen\@clnwd
\newdimen\@clnht

\newif\if@negarg

\def\@whilenoop#1{}
\def\@whiledim#1\do #2{\ifdim #1\relax#2\@iwhiledim{#1\relax#2}\fi}
\def\@iwhiledim#1{\ifdim #1\let\@nextwhile=\@iwhiledim
        \else\let\@nextwhile=\@whilenoop\fi\@nextwhile{#1}}

\def\@whileswnoop#1\fi{}
\def\@whilesw#1\fi#2{#1#2\@iwhilesw{#1#2}\fi\fi}
\def\@iwhilesw#1\fi{#1\let\@nextwhile=\@iwhilesw
         \else\let\@nextwhile=\@whileswnoop\fi\@nextwhile{#1}\fi}

\def\thinlines{\let\@linefnt\tenln \let\@circlefnt\tencirc
  \@wholewidth\fontdimen8\tenln \@halfwidth .5\@wholewidth}
\def\thicklines{\let\@linefnt\tenlnw \let\@circlefnt\tencircw
  \@wholewidth\fontdimen8\tenlnw \@halfwidth .5\@wholewidth}
\thinlines
%%%%%%%%%%%%%%%%%%%%%%%%%%%%%%%%%%%%%%%%%%%%%%%%%%%%%%%%%%%

\PfadD@cke1pt \Pfadd@cke0.5pt
\def\PfadDicke#1{\PfadD@cke#1 \divide\PfadD@cke by2 \Pfadd@cke\PfadD@cke \multiply\PfadD@cke by2}
\long\def\LOOP#1\REPEAT{\def\BODY{#1}\ITERATE}
\def\ITERATE{\BODY \let\next\ITERATE \else\let\next\relax\fi \next}
\let\REPEAT=\fi
\def\Punkt{\hbox{\raise-2pt\hbox to0pt{\hss$\ssize\bullet$\hss}}}
\def\DuennPunkt(#1,#2){\unskip
  \raise#2 \Einheit\hbox to0pt{\hskip#1 \Einheit
          \raise-2.5pt\hbox to0pt{\hss$\bullet$\hss}\hss}}
\def\NormalPunkt(#1,#2){\unskip
  \raise#2 \Einheit\hbox to0pt{\hskip#1 \Einheit
          \raise-3pt\hbox to0pt{\hss\twelvepoint$\bullet$\hss}\hss}}
\def\DickPunkt(#1,#2){\unskip
  \raise#2 \Einheit\hbox to0pt{\hskip#1 \Einheit
          \raise-4pt\hbox to0pt{\hss\fourteenpoint$\bullet$\hss}\hss}}
\def\Kreis(#1,#2){\unskip
  \raise#2 \Einheit\hbox to0pt{\hskip#1 \Einheit
          \raise-4pt\hbox to0pt{\hss\fourteenpoint$\circ$\hss}\hss}}

%%%%%%%%%%%%%%%%%%%%%
%LaTeX line macros%
%%%%%%%%%%%%%%%%%%%%%
\def\Line@(#1,#2)#3{\@xarg #1\relax \@yarg #2\relax
\@linelen=#3\Einheit
\ifnum\@xarg =0 \@vline
  \else \ifnum\@yarg =0 \@hline \else \@sline\fi
\fi}

\def\@sline{\ifnum\@xarg< 0 \@negargtrue \@xarg -\@xarg \@yyarg -\@yarg
  \else \@negargfalse \@yyarg \@yarg \fi
\ifnum \@yyarg >0 \@tempcnta\@yyarg \else \@tempcnta -\@yyarg \fi
\ifnum\@tempcnta>6 \@badlinearg\@tempcnta0 \fi
\ifnum\@xarg>6 \@badlinearg\@xarg 1 \fi
\setbox\@linechar\hbox{\@linefnt\@getlinechar(\@xarg,\@yyarg)}%
\ifnum \@yarg >0 \let\@upordown\raise \@clnht\z@
   \else\let\@upordown\lower \@clnht \ht\@linechar\fi
\@clnwd=\wd\@linechar
\if@negarg \hskip -\wd\@linechar \def\@tempa{\hskip -2\wd\@linechar}\else
     \let\@tempa\relax \fi
\@whiledim \@clnwd <\@linelen \do
  {\@upordown\@clnht\copy\@linechar
   \@tempa
   \advance\@clnht \ht\@linechar
   \advance\@clnwd \wd\@linechar}%
\advance\@clnht -\ht\@linechar
\advance\@clnwd -\wd\@linechar
\@tempdima\@linelen\advance\@tempdima -\@clnwd
\@tempdimb\@tempdima\advance\@tempdimb -\wd\@linechar
\if@negarg \hskip -\@tempdimb \else \hskip \@tempdimb \fi
\multiply\@tempdima \@m
\@tempcnta \@tempdima \@tempdima \wd\@linechar \divide\@tempcnta \@tempdima
\@tempdima \ht\@linechar \multiply\@tempdima \@tempcnta
\divide\@tempdima \@m
\advance\@clnht \@tempdima
\ifdim \@linelen <\wd\@linechar
   \hskip \wd\@linechar
  \else\@upordown\@clnht\copy\@linechar\fi}

\def\@hline{\ifnum \@xarg <0 \hskip -\@linelen \fi
\vrule height\Pfadd@cke width \@linelen depth\Pfadd@cke
\ifnum \@xarg <0 \hskip -\@linelen \fi}

\def\@getlinechar(#1,#2){\@tempcnta#1\relax\multiply\@tempcnta 8
\advance\@tempcnta -9 \ifnum #2>0 \advance\@tempcnta #2\relax\else
\advance\@tempcnta -#2\relax\advance\@tempcnta 64 \fi
\char\@tempcnta}

\def\Vektor(#1,#2)#3(#4,#5){\unskip\leavevmode
  \xcoord#4\relax \ycoord#5\relax
      \raise\ycoord \Einheit\hbox to0pt{\hskip\xcoord \Einheit
         \Vector@(#1,#2){#3}\hss}}

\def\Vector@(#1,#2)#3{\@xarg #1\relax \@yarg #2\relax
\@tempcnta \ifnum\@xarg<0 -\@xarg\else\@xarg\fi
\ifnum\@tempcnta<5\relax
\@linelen=#3\Einheit
\ifnum\@xarg =0 \@vvector
  \else \ifnum\@yarg =0 \@hvector \else \@svector\fi
\fi
\else\@badlinearg\fi}

\def\@hvector{\@hline\hbox to 0pt{\@linefnt
\ifnum \@xarg <0 \@getlarrow(1,0)\hss\else
    \hss\@getrarrow(1,0)\fi}}

\def\@vvector{\ifnum \@yarg <0 \@downvector \else \@upvector \fi}

\def\@svector{\@sline
\@tempcnta\@yarg \ifnum\@tempcnta <0 \@tempcnta=-\@tempcnta\fi
\ifnum\@tempcnta <5
  \hskip -\wd\@linechar
  \@upordown\@clnht \hbox{\@linefnt  \if@negarg
  \@getlarrow(\@xarg,\@yyarg) \else \@getrarrow(\@xarg,\@yyarg) \fi}%
\else\@badlinearg\fi}

\def\@upline{\hbox to \z@{\hskip -.5\Pfadd@cke \vrule width \Pfadd@cke
   height \@linelen depth \z@\hss}}

\def\@downline{\hbox to \z@{\hskip -.5\Pfadd@cke \vrule width \Pfadd@cke
   height \z@ depth \@linelen \hss}}

\def\@upvector{\@upline\setbox\@tempboxa\hbox{\@linefnt\char'66}\raise
     \@linelen \hbox to\z@{\lower \ht\@tempboxa\box\@tempboxa\hss}}

\def\@downvector{\@downline\lower \@linelen
      \hbox to \z@{\@linefnt\char'77\hss}}

\def\@getlarrow(#1,#2){\ifnum #2 =\z@ \@tempcnta='33\else
\@tempcnta=#1\relax\multiply\@tempcnta \sixt@@n \advance\@tempcnta
-9 \@tempcntb=#2\relax\multiply\@tempcntb \tw@
\ifnum \@tempcntb >0 \advance\@tempcnta \@tempcntb\relax
\else\advance\@tempcnta -\@tempcntb\advance\@tempcnta 64
\fi\fi\char\@tempcnta}

\def\@getrarrow(#1,#2){\@tempcntb=#2\relax
\ifnum\@tempcntb < 0 \@tempcntb=-\@tempcntb\relax\fi
\ifcase \@tempcntb\relax \@tempcnta='55 \or
\ifnum #1<3 \@tempcnta=#1\relax\multiply\@tempcnta
24 \advance\@tempcnta -6 \else \ifnum #1=3 \@tempcnta=49
\else\@tempcnta=58 \fi\fi\or
\ifnum #1<3 \@tempcnta=#1\relax\multiply\@tempcnta
24 \advance\@tempcnta -3 \else \@tempcnta=51\fi\or
\@tempcnta=#1\relax\multiply\@tempcnta
\sixt@@n \advance\@tempcnta -\tw@ \else
\@tempcnta=#1\relax\multiply\@tempcnta
\sixt@@n \advance\@tempcnta 7 \fi\ifnum #2<0 \advance\@tempcnta 64 \fi
\char\@tempcnta}
%%%%%%%%%%%%%%%%%%%%%%%%%%%%%%%%%%%%%%%%%%%%%%%%%%%%%%%%%%%%%

\def\Diagonale(#1,#2)#3{\unskip\leavevmode
  \xcoord#1\relax \ycoord#2\relax
      \raise\ycoord \Einheit\hbox to0pt{\hskip\xcoord \Einheit
         \Line@(1,1){#3}\hss}}
\def\AntiDiagonale(#1,#2)#3{\unskip\leavevmode
  \xcoord#1\relax \ycoord#2\relax %\advance\xcoord by -0.05\relax
      \raise\ycoord \Einheit\hbox to0pt{\hskip\xcoord \Einheit
         \Line@(1,-1){#3}\hss}}
\def\Pfad(#1,#2),#3\endPfad{\unskip\leavevmode
  \xcoord#1 \ycoord#2 \thicklines\ZeichnePfad#3\endPfad\thinlines}
\def\ZeichnePfad#1{\ifx#1\endPfad\let\next\relax
  \else\let\next\ZeichnePfad
    \ifnum#1=1
      \raise\ycoord \Einheit\hbox to0pt{\hskip\xcoord \Einheit
         \vrule height\Pfadd@cke width1 \Einheit depth\Pfadd@cke\hss}%
      \advance\xcoord by 1
    \else\ifnum#1=2
      \raise\ycoord \Einheit\hbox to0pt{\hskip\xcoord \Einheit
        \hbox{\hskip-\PfadD@cke\vrule height1 \Einheit width\PfadD@cke depth0pt}\hss}%
      \advance\ycoord by 1
    \else\ifnum#1=3
      \raise\ycoord \Einheit\hbox to0pt{\hskip\xcoord \Einheit
         \Line@(1,1){1}\hss}
      \advance\xcoord by 1
      \advance\ycoord by 1
    \else\ifnum#1=4
      \raise\ycoord \Einheit\hbox to0pt{\hskip\xcoord \Einheit
         \Line@(1,-1){1}\hss}
      \advance\xcoord by 1
      \advance\ycoord by -1
    \fi\fi\fi\fi
  \fi\next}
\def\hSSchritt{\leavevmode\raise-.4pt\hbox to0pt{\hss.\hss}\hskip.2\Einheit
  \raise-.4pt\hbox to0pt{\hss.\hss}\hskip.2\Einheit
  \raise-.4pt\hbox to0pt{\hss.\hss}\hskip.2\Einheit
  \raise-.4pt\hbox to0pt{\hss.\hss}\hskip.2\Einheit
  \raise-.4pt\hbox to0pt{\hss.\hss}\hskip.2\Einheit}
\def\vSSchritt{\vbox{\baselineskip.2\Einheit\lineskiplimit0pt
\hbox{.}\hbox{.}\hbox{.}\hbox{.}\hbox{.}}}
\def\DSSchritt{\leavevmode\raise-.4pt\hbox to0pt{%
  \hbox to0pt{\hss.\hss}\hskip.2\Einheit
  \raise.2\Einheit\hbox to0pt{\hss.\hss}\hskip.2\Einheit
  \raise.4\Einheit\hbox to0pt{\hss.\hss}\hskip.2\Einheit
  \raise.6\Einheit\hbox to0pt{\hss.\hss}\hskip.2\Einheit
  \raise.8\Einheit\hbox to0pt{\hss.\hss}\hss}}
\def\dSSchritt{\leavevmode\raise-.4pt\hbox to0pt{%
  \hbox to0pt{\hss.\hss}\hskip.2\Einheit
  \raise-.2\Einheit\hbox to0pt{\hss.\hss}\hskip.2\Einheit
  \raise-.4\Einheit\hbox to0pt{\hss.\hss}\hskip.2\Einheit
  \raise-.6\Einheit\hbox to0pt{\hss.\hss}\hskip.2\Einheit
  \raise-.8\Einheit\hbox to0pt{\hss.\hss}\hss}}
\def\SPfad(#1,#2),#3\endSPfad{\unskip\leavevmode
  \xcoord#1 \ycoord#2 \ZeichneSPfad#3\endSPfad}
\def\ZeichneSPfad#1{\ifx#1\endSPfad\let\next\relax
  \else\let\next\ZeichneSPfad
    \ifnum#1=1
      \raise\ycoord \Einheit\hbox to0pt{\hskip\xcoord \Einheit
         \hSSchritt\hss}%
      \advance\xcoord by 1
    \else\ifnum#1=2
      \raise\ycoord \Einheit\hbox to0pt{\hskip\xcoord \Einheit
        \hbox{\hskip-2pt \vSSchritt}\hss}%
      \advance\ycoord by 1
    \else\ifnum#1=3
      \raise\ycoord \Einheit\hbox to0pt{\hskip\xcoord \Einheit
         \DSSchritt\hss}
      \advance\xcoord by 1
      \advance\ycoord by 1
    \else\ifnum#1=4
      \raise\ycoord \Einheit\hbox to0pt{\hskip\xcoord \Einheit
         \dSSchritt\hss}
      \advance\xcoord by 1
      \advance\ycoord by -1
    \fi\fi\fi\fi
  \fi\next}
\def\Koordinatenachsen(#1,#2){\unskip
 \hbox to0pt{\hskip-.5pt\vrule height#2 \Einheit width.5pt depth1 \Einheit}%
 \hbox to0pt{\hskip-1 \Einheit \xcoord#1 \advance\xcoord by1
    \vrule height0.25pt width\xcoord \Einheit depth0.25pt\hss}}
\def\Koordinatenachsen(#1,#2)(#3,#4){\unskip
 \hbox to0pt{\hskip-.5pt \ycoord-#4 \advance\ycoord by1
    \vrule height#2 \Einheit width.5pt depth\ycoord \Einheit}%
 \hbox to0pt{\hskip-1 \Einheit \hskip#3\Einheit 
    \xcoord#1 \advance\xcoord by1 \advance\xcoord by-#3 
    \vrule height0.25pt width\xcoord \Einheit depth0.25pt\hss}}
\def\Gitter(#1,#2){\unskip \xcoord0 \ycoord0 \leavevmode
  \LOOP\ifnum\ycoord<#2
    \loop\ifnum\xcoord<#1
      \raise\ycoord \Einheit\hbox to0pt{\hskip\xcoord \Einheit\Punkt\hss}%
      \advance\xcoord by1
    \repeat
    \xcoord0
    \advance\ycoord by1
  \REPEAT}
\def\Gitter(#1,#2)(#3,#4){\unskip \xcoord#3 \ycoord#4 \leavevmode
  \LOOP\ifnum\ycoord<#2
    \loop\ifnum\xcoord<#1
      \raise\ycoord \Einheit\hbox to0pt{\hskip\xcoord \Einheit\Punkt\hss}%
      \advance\xcoord by1
    \repeat
    \xcoord#3
    \advance\ycoord by1
  \REPEAT}
\def\Label#1#2(#3,#4){\unskip \xdim#3 \Einheit \ydim#4 \Einheit
  \def\lo{\advance\xdim by-.5 \Einheit \advance\ydim by.5 \Einheit}%
  \def\llo{\advance\xdim by-.25cm \advance\ydim by.5 \Einheit}%
  \def\loo{\advance\xdim by-.5 \Einheit \advance\ydim by.25cm}%
  \def\o{\advance\ydim by.25cm}%
  \def\ro{\advance\xdim by.5 \Einheit \advance\ydim by.5 \Einheit}%
  \def\rro{\advance\xdim by.25cm \advance\ydim by.5 \Einheit}%
  \def\roo{\advance\xdim by.5 \Einheit \advance\ydim by.25cm}%
  \def\l{\advance\xdim by-.30cm}%
  \def\r{\advance\xdim by.30cm}%
  \def\lu{\advance\xdim by-.5 \Einheit \advance\ydim by-.6 \Einheit}%
  \def\llu{\advance\xdim by-.25cm \advance\ydim by-.6 \Einheit}%
  \def\luu{\advance\xdim by-.5 \Einheit \advance\ydim by-.30cm}%
  \def\u{\advance\ydim by-.30cm}%
  \def\ru{\advance\xdim by.5 \Einheit \advance\ydim by-.6 \Einheit}%
  \def\rru{\advance\xdim by.25cm \advance\ydim by-.6 \Einheit}%
  \def\ruu{\advance\xdim by.5 \Einheit \advance\ydim by-.30cm}%
  #1\raise\ydim\hbox to0pt{\hskip\xdim
     \vbox to0pt{\vss\hbox to0pt{\hss$#2$\hss}\vss}\hss}%
}
\catcode`\@=13

\catcode`\@=11
\font@\twelverm=cmr10 scaled\magstep1
\font@\twelveit=cmti10 scaled\magstep1
\font@\twelvebf=cmbx10 scaled\magstep1
\font@\twelvei=cmmi10 scaled\magstep1
\font@\twelvesy=cmsy10 scaled\magstep1
\font@\twelveex=cmex10 scaled\magstep1

\newtoks\twelvepoint@
\def\twelvepoint{\normalbaselineskip15\p@
 \abovedisplayskip15\p@ plus3.6\p@ minus10.8\p@
 \belowdisplayskip\abovedisplayskip
 \abovedisplayshortskip\z@ plus3.6\p@
 \belowdisplayshortskip8.4\p@ plus3.6\p@ minus4.8\p@
 \textonlyfont@\rm\twelverm \textonlyfont@\it\twelveit
 \textonlyfont@\sl\twelvesl \textonlyfont@\bf\twelvebf
 \textonlyfont@\smc\twelvesmc \textonlyfont@\tt\twelvett
%Erg\"anzung des fetten Small-Capitals-Fonts:
%
 \ifsyntax@ \def\big##1{{\hbox{$\left##1\right.$}}}%
  \let\Big\big \let\bigg\big \let\Bigg\big
 \else
  \textfont\z@=\twelverm  \scriptfont\z@=\tenrm  \scriptscriptfont\z@=\sevenrm
  \textfont\@ne=\twelvei  \scriptfont\@ne=\teni  \scriptscriptfont\@ne=\seveni
  \textfont\tw@=\twelvesy \scriptfont\tw@=\tensy \scriptscriptfont\tw@=\sevensy
  \textfont\thr@@=\twelveex \scriptfont\thr@@=\tenex
        \scriptscriptfont\thr@@=\tenex
  \textfont\itfam=\twelveit \scriptfont\itfam=\tenit
        \scriptscriptfont\itfam=\tenit
  \textfont\bffam=\twelvebf \scriptfont\bffam=\tenbf
        \scriptscriptfont\bffam=\sevenbf
  \setbox\strutbox\hbox{\vrule height10.2\p@ depth4.2\p@ width\z@}%
  \setbox\strutbox@\hbox{\lower.6\normallineskiplimit\vbox{%
        \kern-\normallineskiplimit\copy\strutbox}}%
 \setbox\z@\vbox{\hbox{$($}\kern\z@}\bigsize@=1.4\ht\z@
 \fi
 \normalbaselines\rm\ex@.2326ex\jot3.6\ex@\the\twelvepoint@}

\font@\fourteenrm=cmr10 scaled\magstep2
\font@\fourteenit=cmti10 scaled\magstep2
\font@\fourteensl=cmsl10 scaled\magstep2
\font@\fourteensmc=cmcsc10 scaled\magstep2
\font@\fourteentt=cmtt10 scaled\magstep2
\font@\fourteenbf=cmbx10 scaled\magstep2
\font@\fourteeni=cmmi10 scaled\magstep2
\font@\fourteensy=cmsy10 scaled\magstep2
\font@\fourteenex=cmex10 scaled\magstep2
\font@\fourteenmsa=msam10 scaled\magstep2
\font@\fourteeneufm=eufm10 scaled\magstep2
\font@\fourteenmsb=msbm10 scaled\magstep2
\newtoks\fourteenpoint@
\def\fourteenpoint{\normalbaselineskip15\p@
 \abovedisplayskip18\p@ plus4.3\p@ minus12.9\p@
 \belowdisplayskip\abovedisplayskip
 \abovedisplayshortskip\z@ plus4.3\p@
 \belowdisplayshortskip10.1\p@ plus4.3\p@ minus5.8\p@
 \textonlyfont@\rm\fourteenrm \textonlyfont@\it\fourteenit
 \textonlyfont@\sl\fourteensl \textonlyfont@\bf\fourteenbf
 \textonlyfont@\smc\fourteensmc \textonlyfont@\tt\fourteentt
%Erg\"anzung des fetten Small-Capitals-Fonts:
%
 \ifsyntax@ \def\big##1{{\hbox{$\left##1\right.$}}}%
  \let\Big\big \let\bigg\big \let\Bigg\big
 \else
  \textfont\z@=\fourteenrm  \scriptfont\z@=\twelverm  \scriptscriptfont\z@=\tenrm
  \textfont\@ne=\fourteeni  \scriptfont\@ne=\twelvei  \scriptscriptfont\@ne=\teni
  \textfont\tw@=\fourteensy \scriptfont\tw@=\twelvesy \scriptscriptfont\tw@=\tensy
  \textfont\thr@@=\fourteenex \scriptfont\thr@@=\twelveex
        \scriptscriptfont\thr@@=\twelveex
  \textfont\itfam=\fourteenit \scriptfont\itfam=\twelveit
        \scriptscriptfont\itfam=\twelveit
  \textfont\bffam=\fourteenbf \scriptfont\bffam=\twelvebf
        \scriptscriptfont\bffam=\tenbf
  \setbox\strutbox\hbox{\vrule height12.2\p@ depth5\p@ width\z@}%
  \setbox\strutbox@\hbox{\lower.72\normallineskiplimit\vbox{%
        \kern-\normallineskiplimit\copy\strutbox}}%
 \setbox\z@\vbox{\hbox{$($}\kern\z@}\bigsize@=1.7\ht\z@
 \fi
 \normalbaselines\rm\ex@.2326ex\jot4.3\ex@\the\fourteenpoint@}

\font@\seventeenrm=cmr10 scaled\magstep3
\font@\seventeenit=cmti10 scaled\magstep3
\font@\seventeensl=cmsl10 scaled\magstep3
\font@\seventeensmc=cmcsc10 scaled\magstep3
\font@\seventeentt=cmtt10 scaled\magstep3
\font@\seventeenbf=cmbx10 scaled\magstep3
\font@\seventeeni=cmmi10 scaled\magstep3
\font@\seventeensy=cmsy10 scaled\magstep3
\font@\seventeenex=cmex10 scaled\magstep3
\font@\seventeenmsa=msam10 scaled\magstep3
\font@\seventeeneufm=eufm10 scaled\magstep3
\font@\seventeenmsb=msbm10 scaled\magstep3
\newtoks\seventeenpoint@
\def\seventeenpoint{\normalbaselineskip18\p@
 \abovedisplayskip21.6\p@ plus5.2\p@ minus15.4\p@
 \belowdisplayskip\abovedisplayskip
 \abovedisplayshortskip\z@ plus5.2\p@
 \belowdisplayshortskip12.1\p@ plus5.2\p@ minus7\p@
 \textonlyfont@\rm\seventeenrm \textonlyfont@\it\seventeenit
 \textonlyfont@\sl\seventeensl \textonlyfont@\bf\seventeenbf
 \textonlyfont@\smc\seventeensmc \textonlyfont@\tt\seventeentt
%Erg\"anzung des fetten Small-Capitals-Fonts:
%
 \ifsyntax@ \def\big##1{{\hbox{$\left##1\right.$}}}%
  \let\Big\big \let\bigg\big \let\Bigg\big
 \else
  \textfont\z@=\seventeenrm  \scriptfont\z@=\fourteenrm  \scriptscriptfont\z@=\twelverm
  \textfont\@ne=\seventeeni  \scriptfont\@ne=\fourteeni  \scriptscriptfont\@ne=\twelvei
  \textfont\tw@=\seventeensy \scriptfont\tw@=\fourteensy \scriptscriptfont\tw@=\twelvesy
  \textfont\thr@@=\seventeenex \scriptfont\thr@@=\fourteenex
        \scriptscriptfont\thr@@=\fourteenex
  \textfont\itfam=\seventeenit \scriptfont\itfam=\fourteenit
        \scriptscriptfont\itfam=\fourteenit
  \textfont\bffam=\seventeenbf \scriptfont\bffam=\fourteenbf
        \scriptscriptfont\bffam=\twelvebf
  \setbox\strutbox\hbox{\vrule height14.6\p@ depth6\p@ width\z@}%
  \setbox\strutbox@\hbox{\lower.86\normallineskiplimit\vbox{%
        \kern-\normallineskiplimit\copy\strutbox}}%
 \setbox\z@\vbox{\hbox{$($}\kern\z@}\bigsize@=2\ht\z@
 \fi
 \normalbaselines\rm\ex@.2326ex\jot5.2\ex@\the\seventeenpoint@}

\catcode`\@=13

\topmatter 
\title Bijections for hook pair identities
\endtitle 
\author C.~Krattenthaler$^\dagger$
\endauthor 
\affil 
Institut f\"ur Mathematik der Universit\"at Wien,\\
Strudlhofgasse 4, A-1090 Wien, Austria.\\
e-mail: KRATT\@Ap.Univie.Ac.At\\
WWW: \tt http://radon.mat.univie.ac.at/People/kratt
\endaffil
\address Institut f\"ur Mathematik der Universit\"at Wien,
Strudlhofgasse 4, A-1090 Wien, Austria.
\endaddress
%\email KRATT\@Ap.Univie.Ac.At\\
%WWW: \tt http://radon.mat.univie.ac.at/People/kratt\endemail
%\dedicatory \enddedicatory
%\date \enddate
\thanks{$^\dagger$ Research partially supported by the Austrian
Science Foundation FWF, grant P13190-MAT}
\endthanks
\abstract 
Short, bijective proofs of identities for multisets of
`hook pairs' (arm-leg pairs) of the cells of certain diagrams are
given. These hook pair identities were originally found by Regev.
\endabstract
\endtopmatter
\document

\subhead 1. Introduction \endsubhead
In their work \cite{\ReVeAA} on asymptotic analysis of degrees of
sequences of symmetric group characters,
Regev and Vershik obtained some hook formulas, which led them to
conjecture surprising identities for multisets of hooks. These
identities were shortly thereafter proved independently by Bessenrodt
\cite{\BessAC}, Janson \cite{\JanSAA}, and Regev and Zeilberger
\cite{\ReZeAA}. Another such identity was added by Postnikov and
Regev \cite{\PoReAA}. Moving one step ahead, Regev \cite{\RegeAD}
observed that in fact all these identities are not only true as
identities for multisets of hooks, but even as identities for
multisets of the corresponding arm-leg pairs. He called the latter
(and we follow this convention) ``{\it hook pairs}." As is shown in
\cite{\RegeAD},
these identities imply several nice formulas for special evaluations
of Schur and Jack polynomials.
All the aforementioned
identities feature hooks and arm-leg pairs of regions which are
built out of (nonshifted) Ferrers diagrams. Finally, in
\cite{\RegeAE}, Regev provided similar identities for regions resulting
from shifted diagrams. 

Regev proves his multiset identities in \cite{\RegeAD, \RegeAE} by inductive
arguments. (The proofs in \cite{\JanSAA, \PoReAA, \ReZeAA} are also
inductive, only Bessenrodt's argument in \cite{\BessAC} is
combinatorial.) The purpose of this paper is to provide short,
bijective proofs of all these identities. In fact, what
I am going to demonstrate is that there is just one ``master bijection"
out of which all the identities result straightforwardly.

In the next section we provide all the relevant definitions and 
formulate, in Theorems~1--3, three key identities from
\cite{\RegeAD, \RegeAE}, which straightforwardly imply all other multiset 
identities in these two papers (and, thus, all the multiset
identities in \cite{\BessAC, \JanSAA, \PoReAA, \ReVeAA, \ReZeAA}).
In the subsequent Section~3, we present our ``master
bijection," which immediately implies the first of these identities,
Theorem~1. The resulting proofs of Theorems~2 and 3 are
then given in Sections~4 and 5.

\subhead 2. Identities for multisets of hook pairs\endsubhead
We recall some basic partition terminology (cf\. \cite{\MacdAC,
Ch.~I, Sec.~1}).
A {\it partition} is a sequence $\mu=(\mu_1,\mu_2,\dots,\mu_\ell)$ of
positive integers, arranged in weakly decreasing order. We identify a
partition $\mu$ with its {\it Ferrers diagram}, which is an
array of cells with $\ell$ left-justified rows and $\mu_i$
cells in row $i$. To each cell, $c$ say, we associate two numbers,
the {\it arm length} and the {\it leg length} of $c$. The arm length
$a(c)$ of $c$
is the number of cells which are (strictly) to the right of $c$ and 
in the same row as $c$. Similarly, the leg length $l(c)$ of $c$
is the number of cells which are (strictly) below of $c$ and 
in the same column as $c$. We call the arm-leg pair $(a(c),l(c))$ of
a cell $c$ the {\it hook pair} of $c$.
For example, Figure~1 shows the Ferrers
diagram $(5,3,3,1)$. The marked cell has arm length $3$ and leg
length $2$, thus, the corresponding hook pair being $(3,2)$. We adopt
the convention of writing $a^b$ for a part $a$ of a partition which occurs $b$
times, so that the partition $(5,3,3,1)$ could also be written as
$(5,3^2,1)$.

\midinsert
\vskip10pt
\vbox{
$$
\smatrix \format \sa\c\quad \s\c\quad \s\c\quad \s\c\quad \s\c\quad \se\\
\hlinefor{11}\\
&\vphantom{f}&& \hbox to0pt{$\times$\hss}&& && && &\\
\hlinefor{11}\\
&\vphantom{f}&& && &\\
\hlinefor7\\
&\vphantom{f}&& && &\\
\hlinefor7\\
&\vphantom{f} &\\
\hlinefor3
\endsmatrix
$$
\centerline{\eightpoint Figure 1}
}
\vskip10pt
\endinsert

Now let the partition $\mu=(\mu_1,\mu_2,\dots,\mu_\ell)$ be given,
and choose $k$ and $n$ such that $\mu$ fits in the $n\times k$
rectangular partition $(k^n)$, i.e., $k\ge \mu_1$ and $n\ge \ell$.
Let us denote the rectangle $(k^n)$ by $R_{n,k}$.
Under these assumptions, we form two skew
diagrams, $\SR_{n,k}(\mu)$ and $\SR_{n,k}(\widetilde\mu)$. The skew
diagram
$\SR_{n,k}(\mu)$ is formed by starting with the $k\times n$ rectangle
$R_{n,k}$, removing a copy of the partition $\mu$ from the top-left
corner of $R_{n,k}$, and then gluing a copy of $\mu$ to the right of
$R_{n,k}$ so that the first row of this copy of $\mu$ is glued to
the first row of $R_{n,k}$. See Figure~2.b for an example with
$k=6$, $n=4$, $\mu=(5,2,1)$. The skew diagram 
$\SR_{n,k}(\widetilde\mu)$ is formed in a similar way. Again one
starts with the $k\times n$ rectangle
$R_{n,k}$. But now a copy of the partition $\mu$ {\it which is rotated by
$180^\circ$} is removed from the {\it bottom-right\/}
corner of $R_{n,k}$, and then a copy of $\mu$ {\it which is rotated by
$180^\circ$} is glued {\it to the left\/} of
$R_{n,k}$ so that the last row of this rotated copy of $\mu$ is glued to
the last row of $R_{n,k}$. See Figure~2.c for an example with
$k=6$, $n=4$, $\mu=(5,2,1)$.

\midinsert
\vskip10pt
\vbox{
$$
\Einheit.6cm
\PfadDicke{1.2pt}
\Pfad(0,0),1111112222\endPfad
\Pfad(0,0),2222111111\endPfad
\PfadDicke{.5pt}
\Pfad(0,1),111111\endPfad
\Pfad(0,2),111111\endPfad
\Pfad(0,3),111111\endPfad
\Pfad(1,0),2222\endPfad
\Pfad(2,0),2222\endPfad
\Pfad(3,0),2222\endPfad
\Pfad(4,0),2222\endPfad
\Pfad(5,0),2222\endPfad
\hbox{\hskip5cm}
\PfadDicke{1.2pt}
\Pfad(0,0),1111112222\endPfad
\Pfad(0,0),2\endPfad
\PfadDicke{.5pt}
\SPfad(0,1),21\endSPfad
\SPfad(0,2),211\endSPfad
\SPfad(0,3),211111\endSPfad
\SPfad(1,2),22\endSPfad
\SPfad(2,3),2\endSPfad
\SPfad(3,3),2\endSPfad
\SPfad(4,3),2\endSPfad
\Pfad(0,1),11111112\endPfad
\Pfad(1,2),11111112\endPfad
\Pfad(2,3),111111\endPfad
\Pfad(1,0),22\endPfad
\Pfad(2,0),2221111111112\endPfad
\Pfad(3,0),222\endPfad
\Pfad(4,0),222\endPfad
\Pfad(5,0),2222111111\endPfad
\Pfad(6,1),222\endPfad
\Pfad(7,2),22\endPfad
\Pfad(8,3),2\endPfad
\Pfad(9,3),2\endPfad
\Pfad(10,3),2\endPfad
\hskip6.6cm
$$
\centerline{\eightpoint a. The rectangle $R_{4,6}$\hskip2cm
b. The skew diagram $\SR_{4,6}((5,2,1))$\hskip1cm}
$$
\Einheit.6cm
\PfadDicke{1.2pt}
\Pfad(0,0),2222111111\endPfad
\Pfad(0,0),1\endPfad
\Pfad(6,3),2\endPfad
\PfadDicke{.5pt}
\SPfad(1,0),11111222\endSPfad
\SPfad(2,0),2\endSPfad
\SPfad(3,0),2\endSPfad
\SPfad(4,0),211\endSPfad
\SPfad(5,0),221\endSPfad
\Pfad(-5,0),11111\endPfad
\Pfad(-5,0),2111111111222\endPfad
\Pfad(-4,0),2\endPfad
\Pfad(-3,0),2\endPfad
\Pfad(-2,0),22111111122\endPfad
\Pfad(-1,0),2221111111\endPfad
\Pfad(1,0),2222\endPfad
\Pfad(2,1),222\endPfad
\Pfad(3,1),222\endPfad
\Pfad(4,1),222\endPfad
\Pfad(5,2),22\endPfad
\Pfad(6,3),2\endPfad
\hskip1cm
$$
\centerline{\eightpoint c. The skew diagram
$\SR_{4,6}(\widetilde{(5,2,1)})$}
\vskip7pt
\centerline{\eightpoint Figure 2}
}
\vskip10pt
\endinsert

Now we are ready to state the first of Regev's hook pair identities
\cite{\RegeAD, Theorem~2}.
\proclaim{Theorem 1}The multiset of hook pairs of the cells of
$\SR_{n,k}(\mu)$ is identical with the multiset of hook pairs of the cells of
$\SR_{n,k}(\widetilde\mu)$.
\endproclaim

\midinsert
\vskip10pt
\vbox{
$$
\Einheit.6cm
\PfadDicke{1.2pt}
\Pfad(0,0),2222111111\endPfad
\Pfad(0,0),11\endPfad
\Pfad(6,3),2\endPfad
\PfadDicke{.5pt}
\SPfad(2,0),1111222\endSPfad
\SPfad(3,0),2\endSPfad
\SPfad(4,0),211\endSPfad
\SPfad(5,0),221\endSPfad
\Pfad(-4,0),111112\endPfad
\Pfad(-4,0),211111111222\endPfad
\Pfad(-3,0),2\endPfad
\Pfad(-2,0),22111111122\endPfad
\Pfad(-1,0),2221111111\endPfad
\Pfad(1,0),2222\endPfad
\Pfad(2,0),222221111\endPfad
\Pfad(3,1),2222\endPfad
\Pfad(4,1),2222211\endPfad
\Pfad(5,2),222221\endPfad
\Pfad(6,3),2222\endPfad
\hskip1cm
$$
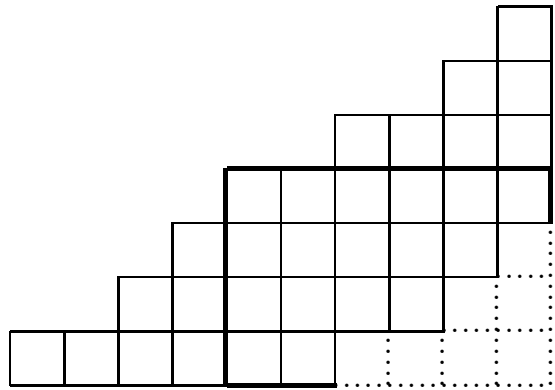
\centerline{\eightpoint The skew diagram
$\SQ(4,6,(4,2,1))$}
\vskip7pt
\centerline{\eightpoint Figure 3}
}
\vskip10pt
\endinsert

Still given a partition $\mu=(\mu_1,\mu_2,\dots,\mu_\ell)$ and
integers $k$ and $n$ such that $k\ge\mu_1$ and $n\ge\ell$, we form
another skew diagram, $\SQ(n,k,\mu)$, by starting again with the
$n\times k$ rectangle $R_{n,k}$, removing a copy of $\mu$ which is
rotated by $180^\circ$ from the bottom-right corner of $R_{n,k}$,
gluing such a rotated copy of $\mu$ to the left of
$R_{n,k}$ in the same way as before when we formed
$\SR_{n,k}(\widetilde\mu)$, and finally gluing a rotated (by
$180^\circ$) copy of $\mu$ to the top of $R_{n,k}$ so that the last
column of this rotated copy of $\mu$ is glued to the last column of
$R_{n,k}$. See Figure~3 for an example with 
$k=6$, $n=4$, $\mu=(4,2,1)$.

The second of Regev's hook pair identities \cite{\RegeAD,
Theorem~1.(a)} reads as follows.

\proclaim{Theorem 2}The multiset of hook pairs of the cells of
$\SQ(n,k,\mu)$ is equal to the union of 
the multiset of hook pairs of the cells of
$R_{n,k}$ and the multiset of hook pairs of the cells of $\mu$.
\endproclaim

\midinsert
\vskip10pt
\vbox{
$$
\Einheit.2cm
\PfadDicke{1.5pt}
\Pfad(0,24),2222222222222222222221\endPfad
\Pfad(0,24),111222222111211111122\endPfad
\Pfad(1,44),2\endPfad
\Pfad(1,44),1\endPfad
\Pfad(2,43),2\endPfad
\Pfad(2,43),1\endPfad
\Pfad(3,42),2\endPfad
\Pfad(3,42),1\endPfad
\Pfad(4,41),2\endPfad
\Pfad(4,41),1\endPfad
\Pfad(5,40),2\endPfad
\Pfad(5,40),1\endPfad
\Pfad(6,39),2\endPfad
\Pfad(6,39),1\endPfad
\Pfad(7,38),2\endPfad
\Pfad(7,38),1\endPfad
\Pfad(8,37),2\endPfad
\Pfad(8,37),1\endPfad
\Pfad(9,36),2\endPfad
\Pfad(9,36),1\endPfad
\Pfad(10,35),2\endPfad
\Pfad(10,35),1\endPfad
\Pfad(11,34),2\endPfad
\Pfad(11,34),1\endPfad
\Pfad(12,33),2\endPfad
\PfadDicke{.7pt}
\Pfad(1,45),111111111111111111111\endPfad
\Pfad(12,33),1112222221222111111222\endPfad
\Label\r{p(\mu)}(3,36)
\Label\r{\mu}(15,33)
\PfadDicke{1.5pt}
\Pfad(1,20),2111111111111111111111\endPfad
\Pfad(1,20),1\endPfad
\Pfad(2,19),2\endPfad
\Pfad(2,19),1\endPfad
\Pfad(3,18),2\endPfad
\Pfad(3,18),1\endPfad
\Pfad(4,17),2\endPfad
\Pfad(4,17),1\endPfad
\Pfad(5,16),2\endPfad
\Pfad(5,16),1\endPfad
\Pfad(6,15),2\endPfad
\Pfad(6,15),1\endPfad
\Pfad(7,14),2\endPfad
\Pfad(7,14),1\endPfad
\Pfad(8,13),2\endPfad
\Pfad(8,13),1\endPfad
\Pfad(9,12),2\endPfad
\Pfad(9,12),1\endPfad
\Pfad(10,11),2\endPfad
\Pfad(10,11),1\endPfad
\Pfad(11,10),2\endPfad
\Pfad(11,10),1\endPfad
\Pfad(12,9),2\endPfad
\Pfad(12,9),1\endPfad
\Pfad(13,8),2\endPfad
\Pfad(13,8),1\endPfad
\Pfad(14,7),2\endPfad
\Pfad(14,7),1\endPfad
\Pfad(15,6),2\endPfad
\Pfad(15,6),1\endPfad
\Pfad(16,5),2\endPfad
\Pfad(16,5),1\endPfad
\Pfad(17,4),2\endPfad
\Pfad(17,4),1\endPfad
\Pfad(18,3),2\endPfad
\Pfad(18,3),1\endPfad
\Pfad(19,2),2\endPfad
\Pfad(19,2),1\endPfad
\Pfad(20,1),2\endPfad
\Pfad(20,1),1\endPfad
\Pfad(21,0),2\endPfad
\Pfad(21,0),1222222222222222222222\endPfad
\PfadDicke{.7pt}
\Pfad(0,0),222222222222222222222111\endPfad
\Pfad(0,0),111111111111111111111\endPfad
\Label\r{q(R(a))}(14,14)
\Label\r{R(a)}(16,23)
\hbox{\hskip6cm}
\PfadDicke{1.5pt}
\Pfad(0,21),2221111112221222222111221111112111222222111\endPfad
\Pfad(0,21),1\endPfad
\Pfad(1,20),2\endPfad
\Pfad(1,20),1\endPfad
\Pfad(2,19),2\endPfad
\Pfad(2,19),1\endPfad
\Pfad(3,18),2\endPfad
\Pfad(3,18),1\endPfad
\Pfad(4,17),2\endPfad
\Pfad(4,17),1\endPfad
\Pfad(5,16),2\endPfad
\Pfad(5,16),1\endPfad
\Pfad(6,15),2\endPfad
\Pfad(6,15),1\endPfad
\Pfad(7,14),2\endPfad
\Pfad(7,14),1\endPfad
\Pfad(8,13),2\endPfad
\Pfad(8,13),1\endPfad
\Pfad(9,12),2\endPfad
\Pfad(9,12),1221111112111222222111222222222222222222222\endPfad
\SPfad(1,21),111111111111111111\endSPfad
\Label\ro{q(A)}(11,16)
\Label\r{A_2}(13,27)
\PfadDicke{.7pt}
\Pfad(0,0),1111111111111111111111222222222222222222222\endPfad
\Pfad(0,0),222222222222222222222\endPfad
\Pfad(0,3),111111222122222211\endPfad
\Pfad(10,11),2\endPfad
\Pfad(10,11),1\endPfad
\Pfad(11,10),2\endPfad
\Pfad(11,10),1\endPfad
\Pfad(12,9),2\endPfad
\Pfad(12,9),1\endPfad
\Pfad(13,8),2\endPfad
\Pfad(13,8),1\endPfad
\Pfad(14,7),2\endPfad
\Pfad(14,7),1\endPfad
\Pfad(15,6),2\endPfad
\Pfad(15,6),1\endPfad
\Pfad(16,5),2\endPfad
\Pfad(16,5),1\endPfad
\Pfad(17,4),2\endPfad
\Pfad(17,4),1\endPfad
\Pfad(18,3),2\endPfad
\Pfad(18,3),1\endPfad
\Pfad(19,2),2\endPfad
\Pfad(19,2),1\endPfad
\Pfad(20,1),2\endPfad
\Pfad(20,1),1\endPfad
\Pfad(21,0),2\endPfad
\Label\lo{\SQ(a,\mu)}(6,34)
\hskip4.8cm
$$
\centerline{\eightpoint Figure 4}
}
\vskip10pt
\endinsert

Finally we concern ourselves with Regev's refinement \cite{\RegeAE}
of Theorem~2 for partitions $\mu$ of the form
$\mu=(\la_1,\dots,\la_s\mid \la_1-1,\dots,\la_s-1)$ (in Frobenius
notation; cf\. \cite{\MacdAC, p.~3}), where
$\la_1>\dots>\la_s>0$. Given such a partition $\mu$, we split it into two
``halves" by cutting it along the diagonal as is illustrated in the
top-left part of
Figure~4 in an example where $\mu=(22^3,16^3,15^6,12^2,6,3^6)=
(21,20,19,12,11,10,8,7,6,5,4,3\mid 20,19,18,11,10,9,7,6,5,4,3,2)$. 
Denote the lower-left region by $p(\mu)$. Given
$a\ge\mu_1-1=\la_1$, let $R(a)$ denote the $a\times(a+1)$ rectangle
$R_{a,a+1}=((a+1)^a)$. Again, we split it into two halves, by cutting
it along the diagonal as is illustrated in the bottom-left part of
Figure~4 for $a=21$. Denote the
upper-right region by $q(R(a))$. Now consider
$\SQ(a,\mu):=\SQ(a,a+1,\mu)$. Recall that this region consists of the
$a\times(a+1)$ rectangle $R_{a,a+1}$, of which a rotated copy of $\mu$
has been removed from its bottom-right corner, and on top of which and to
the left of which was placed a rotated copy of $\mu$ each. Once again we
split it into two halves along the diagonal of the rectangle
$R_{a,a+1}$, as is illustrated in Figure~4.
(There, we have chosen $a=21$ and 
$\mu=(22^3,16^3,15^6,12^2,6,3^6)$. Since it is of no relevance for us
here, we have omitted to display the rotated copy of $\mu$ placed to
the left of the rectangle $R_{21,22}$.) The part of $\SQ(a,\mu)$ above
the diagonal (in the example of Figure~4 this is the region inside
the thick boundary) consists of two
parts, the rotated copy of $\mu$ on top, which we denote by $A_2$, and the
part which remained from $R(a)$, which we denote by $q(A)$. 

If $H$ is a subregion of $G$ (an example being $H=p(\mu)$ and
$G=\mu$), let us write $\HP_G(H)$ for the multiset of hook pairs of
the cells of $H$ {\it measured inside $G$}, i.e., arm length and
leg length are taken with respect to the boundaries of $G$ 
(and not with respect to the boundaries of $H$). (See \cite{\RegeAE,
Sec.~1} for an elaborate example. The definition is motivated by the
definition of {\it shifted\/} hook length when
$p(\mu)$ is regarded as a {\it shifted\/} partition, cf\.
\cite{\MacdAC, Ch.~III, Sec.~8, Ex.~12}.)

With this notation, 
the main result from \cite{\RegeAE, Theorem~I, (I.1)} reads as
follows.

\proclaim{Theorem 3}With the assumptions and notations as explained
above, the following identity holds between multisets of
hook pairs:
$$\HP_{\mu}(p(\mu))\cup \HP_{R(a)}(q(R(a)))=
\HP_{\SQ(a,\mu)}(q(A))\cup \HP_{\SQ(a,\mu)}(A_2).$$
\endproclaim

\subhead 3. The ``master bijection" --- Proof of Theorem~1 \endsubhead
It is obvious that for a proof of Theorem~1 it suffices to consider
just those cells $c$ of $\SR_{n,k}(\mu)$ and $\SR_{n,k}(\widetilde\mu)$
for which the corresponding hook pairs $(a(c),l(c))$ have a fixed arm
length, $a(c)=d$ say, and to show that the multiset of leg lengths of
these cells in $\SR_{n,k}(\mu)$ agrees with the multiset of leg lengths of
these cells in $\SR_{n,k}(\widetilde\mu)$. To illustrate what we mean,
choose $n=10$, $k=8$, $\mu=(7,7,5,4,4,3,3,1)$. The skew diagrams 
$\SR_{n,k}(\mu)$ and $\SR_{n,k}(\widetilde\mu)$ with this choice of
parameters are displayed in Figure~5. The Figure also shows the cells
in $\SR_{n,k}(\mu)$ and $\SR_{n,k}(\widetilde\mu)$ which have arm
length $d=2$. The numbers inside the cells are the corresponding leg
lengths. We call the cells of a region $R$ which have arm length $d$, 
the {\it broken column of $R$ in distance $d$}. We have to show that, in
general, for any $d$ the multiset of leg lengths of cells in the
broken column of $\SR_{n,k}(\mu)$ in distance $d$ is equal to 
the multiset of leg lengths of cells in the
broken column of $\SR_{n,k}(\widetilde\mu)$ in distance $d$.

\midinsert
\vskip10pt
\vbox{
$$
\Einheit.4cm
\PfadDicke{1.5pt}
\Pfad(-8,0),1111111122121122122121122\endPfad
\Pfad(-8,0),2212112212212112211111111\endPfad
\PfadDicke{.5pt}
\Pfad(-3,0),22121\endPfad
\Pfad(-2,0),2212\endPfad
\Pfad(-3,1),1\endPfad
\Pfad(0,3),22122121\endPfad
\Pfad(0,3),12212212\endPfad
\Pfad(0,4),1\endPfad
\Pfad(1,6),1\endPfad
\Pfad(4,8),22\endPfad
\Pfad(4,8),122\endPfad
\Pfad(4,9),1\endPfad
\Label\ro{0}(-3,0)
\Label\ro{1}(-3,1)
\Label\ro{2}(-2,2)
\Label\ro{1}(0,3)
\Label\ro{2}(0,4)
\Label\ro{2}(1,5)
\Label\ro{3}(1,6)
\Label\ro{4}(2,7)
\Label\ro{1}(4,8)
\Label\ro{2}(4,9)
\Label\lo{\SR_{n,k}(\mu)}(-6,7)
\hbox{\hskip3cm}
\PfadDicke{1.5pt}
\Pfad(0,0),1111111122112122122112122\endPfad
\Pfad(0,0),2211212212211212211111111\endPfad
\PfadDicke{.5pt}
\Pfad(5,0),221\endPfad
\Pfad(6,0),22\endPfad
\Pfad(5,1),1\endPfad
\Pfad(7,2),21\endPfad
\Pfad(7,2),12\endPfad
\Pfad(8,3),221\endPfad
\Pfad(8,3),122\endPfad
\Pfad(8,4),1\endPfad
\Pfad(9,5),221\endPfad
\Pfad(9,5),122\endPfad
\Pfad(9,6),1\endPfad
\Pfad(11,7),2122\endPfad
\Pfad(11,7),12122\endPfad
\Pfad(12,9),1\endPfad
\Label\ro{0}(5,0)
\Label\ro{1}(5,1)
\Label\ro{2}(7,2)
\Label\ro{1}(8,3)
\Label\ro{2}(8,4)
\Label\ro{3}(9,5)
\Label\ro{4}(9,6)
\Label\ro{2}(11,7)
\Label\ro{1}(12,8)
\Label\ro{2}(12,9)
\Label\lo{\SR_{n,k}(\widetilde\mu)}(5,10)
\hskip3cm
$$
\centerline{\eightpoint Figure 5}
}
\vskip10pt
\endinsert

If we rotate $\SR_{n,k}(\widetilde\mu)$ by $180^\circ$, then 
$\SR_{n,k}(\mu)$ and the rotated 
$\SR_{n,k}(\widetilde\mu)$ fit together side by
side. Figure 6 shows the result in the case of our example of
Figure~5. The numbers in the broken column which is to the left of
the staircase that forms the border between the two regions are
the vertical 
distances of the cells to the ``bottom" of the diagram (consisting of
the staircase and the base line to the left of the staircase), while the
numbers in the broken column which is to the right of
the staircase are the vertical distances of the cells to the 
``top" of the diagram (consisting of
the staircase and the top line to the right of the staircase).
Our task is to set up a bijection between these two multisets of
numbers.

\midinsert
\vskip10pt
\vbox{
$$
\PfadDicke{2pt}
\Pfad(-4,0),111111111111111\endPfad
\Pfad(0,10),1111111111111\endPfad
\Pfad(0,0),22121122122121122\endPfad
\PfadDicke{.6pt}
\Pfad(-3,0),22121\endPfad
\Pfad(-2,0),2212\endPfad
\Pfad(-3,1),1\endPfad
\Pfad(0,3),22122121\endPfad
\Pfad(0,3),12212212\endPfad
\Pfad(0,4),1\endPfad
\Pfad(1,6),1\endPfad
\Pfad(4,8),22\endPfad
\Pfad(4,8),122\endPfad
\Pfad(4,9),1\endPfad
\Pfad(2,0),22121\endPfad
\Pfad(3,0),2212\endPfad
\Pfad(2,1),1\endPfad
\Pfad(5,3),22122121\endPfad
\Pfad(5,3),12212212\endPfad
\Pfad(5,4),1\endPfad
\Pfad(6,6),1\endPfad
\Pfad(9,8),22\endPfad
\Pfad(9,8),122\endPfad
\Pfad(9,9),1\endPfad
\Label\ro{0}(-3,0)
\Label\ro{1}(-3,1)
\Label\ro{2}(-2,2)
\Label\ro{1}(0,3)
\Label\ro{2}(0,4)
\Label\ro{2}(1,5)
\Label\ro{3}(1,6)
\Label\ro{4}(2,7)
\Label\ro{1}(4,8)
\Label\ro{2}(4,9)
\Label\ro{2}(2,0)
\Label\ro{1}(2,1)
\Label\ro{2}(3,2)
\Label\ro{4}(5,3)
\Label\ro{3}(5,4)
\Label\ro{2}(6,5)
\Label\ro{1}(6,6)
\Label\ro{2}(7,7)
\Label\ro{1}(9,8)
\Label\ro{0}(9,9)
\hskip5cm
$$
\centerline{\eightpoint Figure 6}
}
\vskip10pt
\endinsert

\proclaim{Lemma}
Let $S$ be a given staircase (see Figure~6). Let $C_l$ be the
broken column in distance $d$ left of $S$, and let $C_r$ be
the broken column in distance $d$ right of $S$. Then there is an
explicit bijection (the ``master bijection") between the multiset
$\{l(c):c\in C_l\}$ and the multiset 
$\{l'(c):c\in C_r\}$. Here, $l(c)$ denotes the distance of
cell $c$ to the bottom of the diagram, 
and $l'(c)$ denotes the distance of cell $c$ to the top of the diagram.
\endproclaim

\demo{Proof} 
We claim that the following algorithm defines such a bijection:
\roster
\item"(MB1)" Read the numbers in the cells of $C_l$ one after the other 
by considering the cells in the order bottom to top. While reading, 
out of each maximal increasing subsequence of numbers form a stack.
Place the stacks in such a way that numbers of the same size are at
the same height. (In our example in Figure~6, the numbers
$0,1,2,1,2,2,3,4,1,2$ are read in. See the left part of Figure~7 
for the result after dividing this sequence into stacks.)
\item"(MB2)" Move all the numbers of the first stack which are smaller than
the numbers in the second stack to the second stack. Repeat this
procedure with the second and third stack, etc. (Thus, in our example
in the left part of Figure~7 we would move $0$ from the first to the
second stack, $0,1$ from the second to the third, and $0$ from the
third to the fourth stack. The result is displayed in the right part
of Figure~7.)
\item"(MB3)" Read the numbers in the result in the order from the last stack
to the first, and in each stack from bottom to top. (Thus, from the
right part of Figure~7 we would read $0,1,2,1,2,3,4,2,1,2$.)
\endroster
The claim is that the output of this algorithm is the numbers in the
cells of $C_r$ read by considering the cells in the order top to
bottom (which for our running example can be verified in Figure~6). 
Clearly, this would immediately prove the Lemma because it is
obvious how to invert the algorithm.

\midinsert
\vskip10pt
\vbox{
$$\matrix \format \c&\kern20pt\c&\kern20pt\c&\kern20pt\c\\
&&4\\
&&3&\\
2&2&2&2\\
1&1&&1\\
0
\endmatrix
\hskip4cm
\matrix \format \c&\kern20pt\c&\kern20pt\c&\kern20pt\c\\
&&4\\
&&3\\
2&2&2&2\\
1&&1&1\\
&&&0
\endmatrix$$
\centerline{\eightpoint Figure 7}
}
\vskip10pt
\endinsert

In order to prove the claim, we first introduce some notation. 
The staircase $S$ consists of, alternately, vertical and horizontal
pieces. Let the lengths of these pieces be $v_1,h_1,v_2,h_2,\dots,
v_{p-1},h_{p-1},v_p$,
where $v_i$ stands for the length of the $i$-th vertical piece, if
counted from bottom to top, and where $h_i$ stands for the length of
the $i$-th horizontal piece. In the example in Figure~6 we have
$v_1=2$, $h_1=1$,
$v_2=1$, $h_2=2$,
$v_3=2$, $h_3=1$,
$v_4=2$, $h_4=1$,
$v_5=1$, $h_5=2$,
$v_6=2$. 

Furthermore, let the maximal increasing subsequences when reading the
numbers from $C_l$ (from bottom to top) be $0,1,\dots,M_1$;
$m_2,m_2+1,\dots,M_2$; \dots; $m_q,m_q+1,\dots,M_q$, so that the corresponding
diagram in the style of Figure~7 looks like
$$\matrix \format \c&\kern20pt\c&\kern20pt\c&\kern20pt\c&\kern20pt\c\\
 &M_2& & &M_q\\
 &\vdots&M_3& &\vdots\\
M_1&m_2&\vdots& \dots&m_q\\
\vdots& &m_3& \\
0&
\endmatrix$$
Then there exist uniquely determined integers $i_1,\dots,i_q$ and
$j_1,\dots,j_{q}$ with $1\le i_1<\dots<i_q= p$,
$1= j_1<\dots<j_{q}\le p$, and $j_t\le i_t$ for all $t$, such
that
$$\alignat 3
M_1&=v_{j_1}+\dots+v_{i_1}-1,& h_{j_1}+\dots+h_{i_1-1}&\le d \quad\text{and}&
h_{j_1}+\dots+h_{i_1}&> d,\\
m_2&=v_{j_2}+\dots+v_{i_1},& h_{j_2-1}+\dots+h_{i_1}&> d \quad\text{and}&
h_{j_2}+\dots+h_{i_1}&\le d,\\
M_2&=v_{j_2}+\dots+v_{i_2}-1,& h_{j_2}+\dots+h_{i_2-1}&\le d \quad\text{and}&
h_{j_2}+\dots+h_{i_2}&> d,\tag1\\
\multispan6{\hbox to 11cm{\indent\leaderfill}}\\
m_q&=v_{j_{q}}+\dots+v_{i_{q-1}},& h_{j_{q}-1}+\dots+h_{i_{q-1}}&> d \quad\text{and}&
\quad h_{j_{q}}+\dots+h_{i_{q-1}}&\le d,\\
M_q&=v_{j_{q}}+\dots+v_{i_q}-1,& h_{j_{q}}+\dots+h_{i_q-1}&\le d.
\endalignat$$

Similarly, let the maximal increasing subsequences when reading the
numbers from $C_r$ (from top to bottom) be $0,1,\dots,\tilde M_q$;
\dots; $\tilde m_2,\tilde m_2+1,\dots,\tilde M_2$;
$\tilde m_1,\tilde m_1+1,\dots,\tilde M_1$, 
so that the corresponding
diagram in the style of Figure~7 looks like
$$\matrix \format \c&\kern20pt\c&\kern20pt\c&\kern20pt\c&\kern20pt\c\\
 &\tilde M_2& & &\tilde M_q\\
\tilde M_1 &\vdots& &\tilde M_{q-1}&\vdots\\
\vdots&\tilde m_2&\dots&\vdots&\vdots\\
\tilde m_1& & &\tilde m_{q-1}&\vdots \\
& & & &0
\endmatrix$$
Then there exist uniquely determined integers $\tilde j_1,\dots,\tilde 
j_q$ and
$\tilde i_1,\dots,\tilde i_{q}$ with $1\le \tilde j_1<\dots<\tilde j_q= p$,
$1= \tilde i_1<\dots<\tilde i_{q}\le p$, and $\tilde j_t\ge \tilde i_t$ for all $t$, such
that
$$\alignat 3
\tilde M_q&=v_{\tilde j_q}+\dots+v_{\tilde i_q}-1,\\
&&\hbox{\hskip-2cm} h_{\tilde j_q-1}+\dots+h_{\tilde
i_q}&\le d \quad\text{and}&
h_{\tilde j_q-1}+\dots+h_{\tilde i_q-1}&> d,\\
\tilde m_{q-1}&=v_{\tilde j_{q-1}}+\dots+v_{\tilde i_q},\\
&&\hbox{\hskip-2cm} h_{\tilde
j_{q-1}}+\dots+h_{\tilde i_q-1}&> d \quad\text{and}&
h_{\tilde j_{q-1}-1}+\dots+h_{\tilde i_{q}-1}&\le d,\\
\tilde M_{q-1}&=v_{\tilde j_{q-1}}+\dots+v_{\tilde i_{q-1}}-1,
\hbox{\hskip-2cm}& \\
&&\hbox{\hskip-2cm}h_{\tilde j_{q-1}-1}+\dots+h_{\tilde i_{q-1}}&\le d 
\quad\text{and}&\quad 
h_{\tilde j_{q-1}-1}+\dots+h_{\tilde i_{q-1}-1}&> d,\tag2\\
\multispan6{\hbox to 11cm{\indent\leaderfill}}\\
\tilde m_1&=v_{\tilde j_{1}}+\dots+v_{\tilde i_{2}},\\
&&\hbox{\hskip-2cm} 
h_{\tilde j_{1}}+\dots+h_{\tilde i_{2}-1}&> d \quad\text{and}&
\quad h_{\tilde j_{1}-1}+\dots+h_{\tilde i_{2}-1}&\le d,\\
\tilde M_1&=v_{\tilde j_{1}}+\dots+v_{\tilde i_1}-1,&
\multispan3{\quad $h_{\tilde j_{1}-1}+\dots+h_{\tilde i_1}\le
d.$\hfill}
\endalignat$$

The claim is equivalent to the assertion that $\tilde M_t=M_t$ and
$\tilde m_t=m_{t+1}$ for all $t$. In view of (1) and (2), 
this would
immediately follow once we show that $\tilde j_t=i_t$ and $\tilde
i_t=j_t$ for all $t$. In order to do that, it suffices to derive the
inequalities in (2) from those of (1). Indeed, the general form of
the inequalities in (1) is
$$\alignat 2
h_{j_{t+1}-1}+\dots+h_{i_t}&>d\quad \text {and}&\quad 
h_{j_{t+1}}+\dots+h_{i_t}&\le d,\tag3\\
h_{j_t}+\dots+h_{i_t-1}&\le d\quad \text {and}&
h_{j_t}+\dots+h_{i_t}&> d.\tag4
\endalignat$$
In particular, the first inequality in (3) implies
that 
$$h_{j_{t+1}-1}+\dots+h_{i_{t+1}-1}>d,\tag5$$ 
since by assumption we have $i_t\le i_{t+1}-1$. 
Similarly, the first inequality in (4) implies
that 
$$h_{j_{t+1}-1}+\dots+h_{i_{t}-1}\le d,\tag6$$ 
since by assumption we have $j_t\le j_{t+1}-1$. 
Altogether, the inequalities in (3)--(6) cover those of (2) with
$\tilde j_t=i_t$ and $\tilde i_t=j_t$.

This concludes the proof of the Lemma.

\enddemo

\subhead 4. Proof of Theorem~2 \endsubhead
The preceding bijection allows us to construct a bijection for Theorem~2. We have to
set up a bijection between the hook pairs in $R_{n,k}\cup \mu$ and the
hook pairs in $\SQ(n,k,\mu)$. 

To begin with, we identify the hook pairs in a subregion of $R_{n,k}$
with the hook pairs in 
a subregion of $\SQ(n,k,\mu)$. The subregion of $R_{n,k}$, $S_1$
say, consists of the first $k-\mu_1$ columns of $R_{n,k}$ together with a
copy of $\mu$, reflected upside down, on the bottom of $R_{n,k}$ (see
the left part of 
Figure~8 for an example with $n=18$, $k=24$, $\mu=(15^3,6^3,3^6)$). 
The subregion of $\SQ(n,k,\mu)$, $S_2$ say, consists of the (rotated)
copy of $\mu$ on the left of $\SQ(n,k,\mu)$, together
with the next $k-\mu_1$ columns of $\SQ(n,k,\mu)$ (see the right part
of Figure~8).

\midinsert
\vskip10pt
\vbox{
$$
%\Einheit.6cm
\Pfad(0,0),22222211111111\endPfad
\Pfad(0,0),11111111222222\endPfad
\Pfad(3,4),22\endPfad
\SPfad(3,0),2222\endSPfad
\Pfad(3,4),1\endPfad
\Pfad(4,2),22\endPfad
\Pfad(4,2),1\endPfad
\Pfad(5,1),2\endPfad
\Pfad(5,1),111\endPfad
\Label\o{S_1}(3,1)
\Label\o{T_1}(6,3)
\Label\o{R_{n,k}}(4,6)
\hbox{\hskip5cm}
\Pfad(0,0),11111111211121221222222\endPfad
\Pfad(0,0),21112122122111211121221\endPfad
\Pfad(8,1),22222\endPfad
\SPfad(5,0),2222\endSPfad
\Label\o{S_2}(5,2)
\Label\lo{T_2}(11,4)
\Label\o{\SQ(n,k,\mu)}(8,7)
\hskip6.5cm
$$
\centerline{\eightpoint Figure 8}
}
\vskip10pt
\endinsert

It is completely obvious that, row-wise, the hook pairs in these
subregions must be the same. I.e., reading hook pairs from left
to right in the first row of $S_1$ gives exactly the same 
as reading hook pairs from left to right in the first row of $S_2$,
the same being true for the second rows, etc.

Therefore, what remains is to identify the hook pairs in the
remaining regions. Let us denote the complement of $S_1$ in $R_{n,k}$
by $T_1$, and the complement of $S_2$ in $\SQ(n,k,\mu)$ by $T_2$ (see
Figure~8). Then we have to identify the hook pairs in $T_1\cup \mu$ with the
hook pairs in $T_2$.

Let us again consider broken columns in some given distance $d$ from
the right boundaries. 
This is indicated in Figure~9. (In this example, $n=18$,
$k=24$, $\mu=(15^3,6^3,3^6)$, and $d=4$.) The figure also shows ``fake"
parts of broken columns for $\mu$ and $T_2$, i.e., parts of broken
columns which lie outside of the regions. They are shown with dotted
surroundings and should be ignored for the moment. We have to
identify the multiset of leg lengths in the broken columns of $\mu$
and $T_1$ (excluding the ``fake" part) with the multiset of leg
lengths in the broken column of $T_2$ (excluding the ``fake" part).
Let us denote the multiset of leg lengths in the broken column of
$\mu$ by $L(\mu;d)$, and similarly for $T_1$ and $T_2$.

\midinsert
\vskip10pt
\vbox{
$$
%\Einheit.16666cm
\Einheit.2cm
\PfadDicke{1.5pt}
\Pfad(0,21),222222222222111111111111111\endPfad
\Pfad(0,21),111222222111222111111111222\endPfad
\PfadDicke{.7pt}
\Pfad(10,30),222\endPfad
\Pfad(11,30),222\endPfad
\Pfad(1,27),1222\endPfad
\Pfad(1,27),2221\endPfad
\SPfad(-2,21),222222\endSPfad
\SPfad(-2,27),1\endSPfad
\SPfad(-2,21),1222222\endSPfad
\Label\u{\mu}(6,27)
\PfadDicke{1.5pt}
\Pfad(0,12),222222111111111111111\endPfad
\Pfad(15,3),222222222222222\endPfad
\Pfad(0,12),111\endPfad
\Pfad(3,6),222222\endPfad
\Pfad(3,6),111\endPfad
\Pfad(6,3),222\endPfad
\Pfad(6,3),111111111\endPfad
\SPfad(0,0),222222222222\endSPfad
\SPfad(0,0),111111111111111222\endSPfad
\Label\r{T_1}(15,10)
\PfadDicke{.7pt}
\Pfad(10,3),222222222222222\endPfad
\Pfad(11,3),222222222222222\endPfad
\hbox{\hskip7cm}
\SPfad(0,0),222\endSPfad
\SPfad(0,0),111111111111111222222222222\endSPfad
\PfadDicke{1.5pt}
\Pfad(0,3),111111111222111222222111222222222222222222\endPfad
\Pfad(0,3),222222222222222222111111111222111222222111\endPfad
%\SPfad(0,18),111111111111111\endSPfad
\SPfad(0,12),111111111111\endSPfad
\Label\r{T_2}(15,20)
\PfadDicke{.7pt}
\Pfad(10,12),222222222222\endPfad
\Pfad(10,12),1222222222222\endPfad
\Pfad(7,6),2222221\endPfad
\Pfad(7,6),1222222\endPfad
\Pfad(4,3),2221\endPfad
\Pfad(4,3),1222\endPfad
\SPfad(-5,0),1222\endSPfad
\SPfad(-5,0),222\endSPfad
\SPfad(-5,3),1\endSPfad
\hskip2.5cm
$$
\centerline{\eightpoint Figure 9}
}
\vskip10pt
\endinsert

In order to accomplish this identification, 
we observe that the Lemma from Section~3 provides a bijection
between the multiset of leg lengths in the broken column of $\mu$, 
the ``fake" part {\it included}, and the part of the broken column of
$T_2$ which is to the left of the (rotated) copy of $\mu$ that has
been removed from $\SQ(n,k,\mu)$ (indicated by dotted lines in
Figure~9), again
the ``fake" part {\it included}. (In Figure~9, this part of the
broken column is the one below the dotted horizontal line running
through $T_2$,
the ``fake" part {\it included}.) Let us denote the former multiset
by $L_f(\mu;d)$, and the latter by $L_f(T_2;d)$. 

Furthermore, with the notation $[[N]]:=\{0,1,\dots,N-1\}$, we have
$$\align
&L(\mu;d)\cup L(T_1;d)=\Big(L_f(\mu;d)\backslash
[[\ell(\mu)+1-\min\{k:\mu_k\le d\}]]\Big)\\
&\hskip4cm \quad \cup 
\Big([[n]]\backslash[[\max\{k:\mu_k\ge\mu_1-d\}]]\Big)\\
&\quad =\Big(L_f(\mu;d)\cup [[n]]\Big)\Big\backslash 
\Big([[\ell(\mu)+1-\min\{k:\mu_k\le d\}]] \cup
[[\max\{k:\mu_k\ge\mu_1-d\}]]\Big),
\tag7\endalign$$
and
$$\align 
&L(T_2;d)=\Big(L_f(T_2;d)\backslash [[\max\{k:\mu_k+d\ge\mu_1\}]]\Big)\\
&\hskip4cm \quad 
\cup \Big([[n]]\backslash 
[[\ell(\mu)+1-\min\{k:\mu_k\le d\}]]\Big)\\
&\quad\! =\Big(L_f(T_2;d)\cup [[n]]\Big)\Big\backslash
\Big([[\max\{k:\mu_k+d\ge\mu_1\}]]\cup 
[[\ell(\mu)+1-\min\{k:\mu_k\le d\}]]\Big).
\tag8\endalign$$
In view of our previous observation that the multisets $L_f(\mu;d)$ and 
$L_f(T_2;d)$ are in bijection, it is now completely evident that the multisets
$L(\mu;d)\cup L(T_1;d)$ and $L(T_2;d)$ agree.

\remark{Remark} 
The argument thus far
was not completely bijective. But it could easily be made bijective,
for example, using the simplified version of the involution principle
described in \cite{\StanAP, bottom of p.~80ff}. 
There are other possibilities, but
none of these seem to yield naturally defined mappings.
In particular, these would not add any further insight. It is the
``master bijection" from Section~3 and the computations in (7) and (8)
which really explain why Theorem~2 is true.
\endremark

\subhead 5. Proof of Theorem~3 \endsubhead
In a similar manner, Theorem~3 can be proved.
We have to prove that the hook pairs in
$p(\mu)$ and $q(R(a))$ (measured inside $\mu$ and $R(a)$,
respectively) are the same as the hook pairs in $q(A)$ and
$A_2$ (measured inside $\SQ(a,\mu)$), see Figure~4.

Again, it suffices to consider a broken column in some distance $d$ from the
right boundaries in the corresponding figures, and to show that the leg lengths
along the broken columns agree (see Figure~10, where $d=4$; the
``fake" part of the broken column in $p(\mu)$ should be ignored for the
moment), i.e.,
$$L(p(\mu);d)\cup L(q(R(a));d)=L(q(A);d)\cup L(A_2;d),$$
using the same notation as before. 

\topinsert
\vskip10pt
\vbox{
$$
\Einheit.2cm
\PfadDicke{1.5pt}
\Pfad(0,24),2222222222222222222221\endPfad
\Pfad(0,24),111222222111211111122\endPfad
\Pfad(1,44),2\endPfad
\Pfad(1,44),1\endPfad
\Pfad(2,43),2\endPfad
\Pfad(2,43),1\endPfad
\Pfad(3,42),2\endPfad
\Pfad(3,42),1\endPfad
\Pfad(4,41),2\endPfad
\Pfad(4,41),1\endPfad
\Pfad(5,40),2\endPfad
\Pfad(5,40),1\endPfad
\Pfad(6,39),2\endPfad
\Pfad(6,39),1\endPfad
\Pfad(7,38),2\endPfad
\Pfad(7,38),1\endPfad
\Pfad(8,37),2\endPfad
\Pfad(8,37),1\endPfad
\Pfad(9,36),2\endPfad
\Pfad(9,36),1\endPfad
\Pfad(10,35),2\endPfad
\Pfad(10,35),1\endPfad
\Pfad(11,34),2\endPfad
\Pfad(11,34),1\endPfad
\Pfad(12,33),2\endPfad
\PfadDicke{.7pt}
\Pfad(1,45),111111111111111111111\endPfad
\Pfad(12,33),1112222221222111111222\endPfad
\SPfad(-2,24),1222222\endSPfad
\SPfad(-2,24),2222221\endSPfad
\Pfad(1,30),12\endPfad
\Pfad(1,30),21\endPfad
\Pfad(7,31),122\endPfad
\Pfad(7,31),221\endPfad
\Pfad(10,33),122\endPfad
\Pfad(10,33),221\endPfad
\Label\r{p(\mu)}(3,36)
\Kreis(12,33)
\Label\r{\hbox{\, $M$}}(12,32)
\PfadDicke{1.5pt}
\Pfad(1,20),2111111111111111111111\endPfad
\Pfad(1,20),1\endPfad
\Pfad(2,19),2\endPfad
\Pfad(2,19),1\endPfad
\Pfad(3,18),2\endPfad
\Pfad(3,18),1\endPfad
\Pfad(4,17),2\endPfad
\Pfad(4,17),1\endPfad
\Pfad(5,16),2\endPfad
\Pfad(5,16),1\endPfad
\Pfad(6,15),2\endPfad
\Pfad(6,15),1\endPfad
\Pfad(7,14),2\endPfad
\Pfad(7,14),1\endPfad
\Pfad(8,13),2\endPfad
\Pfad(8,13),1\endPfad
\Pfad(9,12),2\endPfad
\Pfad(9,12),1\endPfad
\Pfad(10,11),2\endPfad
\Pfad(10,11),1\endPfad
\Pfad(11,10),2\endPfad
\Pfad(11,10),1\endPfad
\Pfad(12,9),2\endPfad
\Pfad(12,9),1\endPfad
\Pfad(13,8),2\endPfad
\Pfad(13,8),1\endPfad
\Pfad(14,7),2\endPfad
\Pfad(14,7),1\endPfad
\Pfad(15,6),2\endPfad
\Pfad(15,6),1\endPfad
\Pfad(16,5),2\endPfad
\Pfad(16,5),1\endPfad
\Pfad(17,4),2\endPfad
\Pfad(17,4),1\endPfad
\Pfad(18,3),2\endPfad
\Pfad(18,3),1\endPfad
\Pfad(19,2),2\endPfad
\Pfad(19,2),1\endPfad
\Pfad(20,1),2\endPfad
\Pfad(20,1),1\endPfad
\Pfad(21,0),2\endPfad
\Pfad(21,0),1222222222222222222222\endPfad
\PfadDicke{.7pt}
\Pfad(0,0),222222222222222222222111\endPfad
\Pfad(0,0),111111111111111111111\endPfad
\Pfad(17,4),22222222222222222\endPfad
\Pfad(18,4),22222222222222222\endPfad
\Label\r{q(R(a))}(11,16)
\hbox{\hskip6cm}
\PfadDicke{1.5pt}
\Pfad(0,21),2221111112221222222111221111112111222222111\endPfad
\Pfad(0,21),1\endPfad
\Pfad(1,20),2\endPfad
\Pfad(1,20),1\endPfad
\Pfad(2,19),2\endPfad
\Pfad(2,19),1\endPfad
\Pfad(3,18),2\endPfad
\Pfad(3,18),1\endPfad
\Pfad(4,17),2\endPfad
\Pfad(4,17),1\endPfad
\Pfad(5,16),2\endPfad
\Pfad(5,16),1\endPfad
\Pfad(6,15),2\endPfad
\Pfad(6,15),1\endPfad
\Pfad(7,14),2\endPfad
\Pfad(7,14),1\endPfad
\Pfad(8,13),2\endPfad
\Pfad(8,13),1\endPfad
\Pfad(9,12),2\endPfad
\Pfad(9,12),1221111112111222222111222222222222222222222\endPfad
\SPfad(1,21),111111111111111111\endSPfad
\Label\ro{q(A)}(10,18)
\Label\r{A_2}(11,27)
\PfadDicke{.7pt}
\Pfad(0,0),1111111111111111111111222222222222222222222\endPfad
\Pfad(0,0),222222222222222222222\endPfad
\Pfad(0,3),111111222122222211\endPfad
\Pfad(10,11),2\endPfad
\Pfad(10,11),1\endPfad
\Pfad(11,10),2\endPfad
\Pfad(11,10),1\endPfad
\Pfad(12,9),2\endPfad
\Pfad(12,9),1\endPfad
\Pfad(13,8),2\endPfad
\Pfad(13,8),1\endPfad
\Pfad(14,7),2\endPfad
\Pfad(14,7),1\endPfad
\Pfad(15,6),2\endPfad
\Pfad(15,6),1\endPfad
\Pfad(16,5),2\endPfad
\Pfad(16,5),1\endPfad
\Pfad(17,4),2\endPfad
\Pfad(17,4),1\endPfad
\Pfad(18,3),2\endPfad
\Pfad(18,3),1\endPfad
\Pfad(19,2),2\endPfad
\Pfad(19,2),1\endPfad
\Pfad(20,1),2\endPfad
\Pfad(20,1),1\endPfad
\Pfad(21,0),2\endPfad
\Pfad(17,21),222222222222222\endPfad
\Pfad(17,21),1222222222222222\endPfad
\Pfad(14,15),2222221\endPfad
\Pfad(14,15),1222222\endPfad
\Pfad(11,14),21\endPfad
\Pfad(12,14),2\endPfad
\Label\lo{\SQ(a,\mu)}(6,34)
\Kreis(10,12)
\Label\r{\hbox{\ $M'$}}(10,12)
\hskip4.8cm
$$
\centerline{\eightpoint Figure 10}
}
\vskip10pt
\endinsert

Again, the key is the Lemma from Section~3. 
What Figure~11 shows is obtained by the following
construction: In Figure~10, in the part which shows $p(\mu)$, the
extremal point, denoted by $M$, 
on the ``diagonal" is circled. The corresponding point, denoted by
$M'$, in the part which shows $\SQ(a,\mu)$ is also circled.
The top-most part of the broken column in $p(\mu)$ reaches a certain
height, $h$ say, above $M$. (The height $h$ can be any nonnegative
integer.) In Figure~10, we have
$h=2$. Now cut the Figure which represents $\SQ(a,\mu)$ by the
horizontal line which is $h$ units below $M'$. In addition, cut
off $A_2$. This gives Figure~11 (where $h=2$), upon also adding a broken column in
distance $d$ to the right of the path which runs from bottom-left to
top-right, and upon completing a broken column in distance $d$ to the
left of the path.

\topinsert
\vskip10pt
\vbox{
$$
\Einheit.2cm
\PfadDicke{1.5pt}
\Pfad(0,11),1\endPfad
\Pfad(1,10),2\endPfad
\Pfad(1,10),1\endPfad
\Pfad(2,9),2\endPfad
\Pfad(2,9),1\endPfad
\Pfad(3,8),2\endPfad
\Pfad(3,8),1\endPfad
\Pfad(4,7),2\endPfad
\Pfad(4,7),1\endPfad
\Pfad(5,6),2\endPfad
\Pfad(5,6),1\endPfad
\Pfad(6,5),2\endPfad
\Pfad(6,5),1\endPfad
\Pfad(7,4),2\endPfad
\Pfad(7,4),1\endPfad
\Pfad(8,3),2\endPfad
\Pfad(8,3),1\endPfad
\Pfad(9,2),2\endPfad
\Pfad(9,2),1221111112111222222111\endPfad
\PfadDicke{.7pt}
\Pfad(1,11),111111111111111111\endPfad
\Pfad(22,0),22222222222\endPfad
\Pfad(0,0),22222222222\endPfad
\Pfad(0,0),11111112211\endPfad
\Pfad(7,0),1111\endPfad
\Pfad(10,1),2\endPfad
\Pfad(10,1),1\endPfad
\Pfad(11,0),2\endPfad
\Pfad(11,0),11111111111\endPfad
\Pfad(14,5),2222221\endPfad
\Pfad(14,5),1222222\endPfad
\Pfad(11,4),21\endPfad
\Pfad(12,4),2\endPfad
\Pfad(2,0),221\endPfad
\Pfad(3,0),22\endPfad
\Pfad(5,2),221\endPfad
\Pfad(5,2),122\endPfad
\Pfad(11,1),21\endPfad
\Pfad(12,0),22\endPfad
\Pfad(14,2),22\endPfad
\Pfad(14,2),122\endPfad
\Pfad(20,4),21\endPfad
\Pfad(20,4),12\endPfad
\SPfad(23,5),2222221\endSPfad
\SPfad(23,5),1222222\endSPfad
\Kreis(10,2)
\hskip5.2cm
$$
\centerline{\eightpoint Figure 11}
}
\vskip10pt
\endinsert

Let us denote the leg lengths (measured ``to the bottom," i.e., with
respect to the path which runs from bottom-left to top-right) in the
broken column to the left by $L_f(q(A);d)$. These leg lengths contain
all those which appear in $L(q(A);d)$, plus the additional ones in
the added part of the broken column. More precisely, it is not
difficult to see that (at this point the special form of $\mu$ enters)
$$L_f(q(A);d)=L(q(A);d)\cup [[d]].$$
Furthermore, the broken column to the right of the path is in perfect
correspondence with the broken column of $p(\mu)$ (with the ``fake" part
{\it included}; compare Figure~10), which is seen by rotation by
$180^\circ$.
Let us denote the leg lengths (measured ``to the top," i.e., 
with respect to the path which runs from bottom-left to top-right) in the
broken column to the right by $L_f(p(\mu);d)$. Obviously we have
$$L_f(p(\mu);d)=L(p(\mu);d)\cup [[\ell(\mu)+1-\min\{k:\mu_k\le d\}]].$$
The point of the construction in Figure~11 is that 
the Lemma from Section~3 tells us that $L_f(q(A);d)=L_f(p(\mu);d)$.
Now, if we combine everything then we are done: We have
$$\align 
L(p(\mu);d)\cup {}&L(q(R(a));d)\\
&=\Big(L_f(p(\mu);d)\backslash [[\ell(\mu)+1-\min\{k:\mu_k\le d\}]]\Big)
\cup \Big([[a]]\backslash [[d]]\Big)\\
&=\Big(L_f(p(\mu);d)\cup [[a]]\Big)\Big\backslash
 \Big([[\ell(\mu)+1-\min\{k:\mu_k\le d\}]]\cup [[d]]\Big)
\endalign$$
and
$$\align 
L(q(A);d)\cup{}& L(A_2;d)=
\Big(L_f(q(A);d)\backslash [[d]]\Big)\cup 
\Big([[a]]\backslash [[\ell(\mu)+1-\min\{k:\mu_k\le d\}]]\Big)\\
&=\Big(L_f(q(A);d)\cup [[a]] \Big)\Big\backslash
\Big([[d]]\cup [[\ell(\mu)+1-\min\{k:\mu_k\le d\}]]\Big).
\endalign$$
Clearly, in view of $L_f(q(A);d)=L_f(p(\mu);d)$ these two expressions
agree, as desired. Hence, the assertion of Theorem~3 is established.

\remark{Remark} Again, the above argument 
was not completely bijective. Statements analogous to those in the
Remark at the end of Section~4 apply also here.
\endremark

\enddocument